\documentclass{amsart}
\usepackage{amscd} 
\usepackage{amsfonts} 
\usepackage{amssymb} 
\usepackage{latexsym} 
\input{diagrams}
\newcommand{\ncm}{\newcommand}

\def\C{\mathbb{C}\,}

\def\ME{\mathcal{E}}

\newtheorem{theorem}{Theorem}[section]
\newtheorem{prop}[theorem]{Proposition}
\newtheorem{lemma}[theorem]{Lemma}
\newtheorem{cor}[theorem]{Corollary}
\newtheorem{prob}[theorem]{Problem}
\newtheorem{lem&def}[theorem]{Lemma \& Definition}
\newtheorem{definition}[theorem]{Definition}

\ncm{\End}{\mbox{\rm End}\,}

\def\Hom{\mbox{\rm Hom}\,}

\def\Im{\mbox{\rm Im}\,}

\def\id{\mbox{\rm id}}

\def\into{\hookrightarrow}
\def\to{\rightarrow}

\def\o{\otimes}    
\def\bra{\langle}
\def\ket{\rangle}

\ncm{\rarr}[1]{\stackrel{#1}{\longrightarrow}}
\ncm{\larr}[1]{\stackrel{#1}{\longleftarrow}}
\def\cop{\Delta}

\def\eps{\varepsilon}

\def\du1{\hat 1}

\def\-1{_{(-1)}}
\def\0{_{(0)}}
\def\1{_{(1)}}
\def\2{_{(2)}}
\def\3{_{(3)}}

\def\|{\, | \,}

\def\du1{\hat 1}
\def\lact{\triangleright}
\def\ract{\triangleleft}

\begin{document}

\title{Normal Hopf subalgebras, depth two and Galois extensions}
\author{Lars Kadison}
\address{Matematiska Institutionen \\ G{\" o}teborg 
University \\ 
S-412 96 G{\" o}teborg, Sweden} 
\email{lkadison@c2i.net} 
\date{}
\thanks{}
\subjclass{16W30 (13B05, 20L05, 16S40, 81R50)}  
\date{} 

\begin{abstract} 
Let $S$ be the left $R$-bialgebroid of a depth two extension with
centralizer $R$.  
We show that the left endomorphism ring
of depth two extension, not necessarily balanced, is a left $S$-Galois extension of
$A^{\rm op}$. Looking to examples of depth two, 
we establish that a Hopf subalgebra is normal if and only if
it is a Hopf-Galois extension.  We find a class of examples
of the alternative Hopf algebroids in \cite{BS}. 
We also characterize finite weak Hopf-Galois extensions using an alternate Galois 
canonical mapping with several corollaries: that these are depth two and
that surjectivity of the Galois mapping implies its bijectivity.    
\end{abstract} 
\maketitle

\section{Introduction}

\begin{figure}
$$\begin{diagram}
\ME \o_{\rho(A)} \ME &&  \rTo^{\beta}&& S \o_R \ME \\
\dTo^{\cong} && && \uTo_{\cong}    \\
\ME \o_A \ME & &  && \Hom ({}_BA \o_B A, {}_BA) \\
& \SE_{\cong}& & \NE_{\cong} & \\
&& \Hom({}_B\Hom(\ME_A, A_A), {}_BA) &&
\end{diagram}$$
\caption{Galois map $\beta$ factors through various isomorphisms
in Theorem~\ref{th-endogalext}.}
\end{figure}

Bialgebroids are generalized bialgebras over a noncommutative base ring \cite{BW,KS}. As in the theory of bialgebras,
there are associated to a bialgebroid both module and comodule algebras, smash products and Galois
extensions \cite{BW,Bo, Sz, Karl, LK, EN}.
All these constructions are present given a depth two extension $A \| B$,
where the noncommutative base ring is the centralizer $R$ of the ring extension: the extension is depth two if its tensor-square is centrally projective w.r.t.\ the natural $A$-$B$-bimodule as well as
$B$-$A$-bimodule $A$ \cite{KS, Karl, LK}. An extra condition that the natural module $A_B$
is balanced or faithfully flat ensures that $A \| B$ is a right Galois extension w.r.t.\ the $R$-bialgebroid $T := (A \o_B A)^B$
\cite{Karl}.  
In section~2 we show that this condition is not needed for the left endomorphism ring $\ME := \End {}_BA$ to be a left  $S$-Galois
extension of the right multiplications $\rho(A)$.  The proof involves the commutative diagram below.

In section~3 we show that the bialgebroid $T$ of a depth two extension $A$ over a Kanzaki
separable algebra $B$ 
is a Hopf algebroid of the type in \cite[B\"{o}hm-Szlachanyi]{BS}.  The antipode is very
naturally given by a twist of $(A \o_B A)^B$ utilizing a symmetric separability element.  
These provide then further examples of  non-dual Hopf algebroids, in contrast to the dual Hopf
algebroids $S$ and $T$ of a depth two Frobenius extension \cite{BS}.  They are also Hopf algebroids
with no obvious counterparts among Lu's version of Hopf algebroid \cite{Lu}.  

 A depth two complex subalgebra is a generalization of normal subgroup \cite{KK}.
The question was then raised whether depth two Hopf subalgebras
are precisely the normal Hopf subalgebras ($\supseteq$ in \cite{KK}).  
In a very special case, this is true when the notion of
depth two is narrowed to H-separability \cite{Karl}, an exercise in going up and down with ideals as in commutative
algebra. We  study in section~4 the special case of depth two represented by finite Hopf-Galois extensions:
we show that a Hopf-Galois Hopf subalgebra is normal using a certain algebra epimorphism from the Hopf overalgebra
to the Hopf algebra which is coacting Galois, and comparing dimensions of the kernel with the associated Schneider coalgebras. 

A special case of Galois theory for bialgebroids is weak Hopf-Galois theory \cite{BW,CDG,Karl, LK},
(where Hopf-Galois theory is in turn a special case):  for depth two extensions, each type of Galois extension
occurs as we move from any centralizer to separable centralizers to one-dimensional centralizers.  
Conversely, each type of Galois extension, so long as it is finitely generated, is of depth two \cite{KS, Karl, LK}.
In section~5 we complete the proof that a weak Hopf-Galois extension is left  depth two
by studying the alternative Galois mapping $\beta': A \o_B A \to A\o H$
where $\beta'(a \o a') = a\0 a' \o a\1$. As a corollary we find an interesting factorization of the Galois isomorphism
of a weak Hopf algebra over its target subalgebra.  In a second corollary,  a direct proof is given that a surjective
Galois mapping for an  $H$-extension
 is automatically bijective, if $H$ is a finite dimensional weak Hopf algebra.\footnote{My thanks to Tomasz Brzezinski
for pointing out this area of research.}   Finally, it is shown by somewhat
different means than in \cite{BW} that a weak bialgebra in Galois extension of
its target subalgebra has an antipode reconstructible from the Galois mapping.  We provide some evidence
for more generally a
weak bialgebra, which coacts Galois on an algebra over a field, 
 having an antipode, something which is true for bialgebras by a result of
Schauenburg \cite{Sch}.

\section{Depth two and endomorphism ring Galois extensions}

The basic set-up throughout this section is the following. Let
$A \| B$ be a ring extension with centralizer denoted by $R := C_A(B) = A^B$,
bimodule endomorphism ring $S := \End {}_BA_B$ and $B$-central
tensor-square $T := (A \o_B A)^B$. 
  $T$ has a ring structure induced
from $T \cong \End {}_A A \o_B A_A$ given by 
\begin{equation}
\label{eq: tee mult}
tt' = {t'}^1 t^1 \o t^2 {t'}^2, \ \ \ \ \ 1_T = 1 \o 1,
\end{equation}
where $t = t^1 \o t^2 \in T$ uses a Sweedler notation and suppresses
a possible summation over simple tensors.  Let $\lambda: A \into \End A_B$
denote left multiplication and $\rho: A \into \End {}_BA$ denote right multiplication.
Let $\ME$ denote $\End {}_BA$ and
note that $S \subseteq \mathcal{E}$, a subring under the usual
composition of functions.  Note that $\lambda$ restricts to $R \into S$ and $\rho$ restricts
to $R \into S^{\rm op}$.

We have the notion of an arbitrary bimodule
being centrally projective with respect to a canonical bimodule \cite{Karl}: we say that
a bimodule ${}_AM_B$, where $A$ and $B$ are two arbitrary rings,
is \textit{centrally projective w.r.t.\ a bimodule} ${}_AN_B$, if ${}_AM_B$ is
isomorphic to a direct summand of a finite direct sum of $N$
with itself; in symbols,
if ${}_AM_B \oplus * \cong \oplus^n {}_AN_B$.  This covers the usual notion
of centrally projective $A$-$A$-bimodule $P$ where the canonical $A$-$A$-bimodule
is understood to be the natural bimodule $A$ itself. 

We recall the definition of a \textit{depth two} ring extension $A \| B$ as simply  that  its tensor-square
$A \o_B A$  be centrally projective w.r.t.\ the natural $B$-$A$-bimodule $A$ (left D2)
and the natural $A$-$B$-bimodule $A$ (right D2). 
 A very useful characterization of depth two extension is that an extension is D2 if there exist finitely many paired elements (a left D2 quasibase) $\beta_i \in S$,
$t_i \in T$  and finitely many paired elements (right D2 quasibase) $\gamma_j \in S$, $u_j \in T$ such that
\begin{equation}
\label{eq: d2 quasibase}
a \o a' = \sum_i t_i \beta_i(a)a' = \sum_j a \gamma_j(a')u_j 
\end{equation}
for all $a, a' \in A$ \cite[3.7]{KS}: we fix this notation.  
Centrally projective ring extensions, H-separable
extensions and f.g.\
Hopf-Galois extensions are some of the classes of examples of D2 extension. If $A$ and $B$ are the complex
group algebras corresponding to a subgroup $H < G$ of  a finite group, then $A \| B$ is D2 iff $H$ is
a normal subgroup in $G$ \cite{KK}. 
 In a later section of this paper we show the details that finite weak Hopf-Galois extensions are left and
right D2. More generally, Galois extensions for bialgebroids and their comodule algebras are of depth two \cite{LK}. 

Recall that a left bialgebroid 
$R'$-bialgebroid $S'$ is first of all two rings $R'$ and $S'$ with two commuting maps
$\tilde{s}, \tilde{t}: R' \to S'$, a ring homomorphism and anti-homomorphism resp., commuting in the sense  that $\tilde{s}(r) \tilde{t}(r') = \tilde{t}(r') \tilde{s}(r)$ for all $r,r' \in R'$.  Second, it is an $R'$-coring 
$(S',\ \cop \!: S' \to \ $ $ S' \o_{R'} S', \eps: S' \to R')$  \cite{BW} w.r.t.\ the $R'$-$R'$-bimodule $r \cdot x \cdot r' =  \tilde{t}(r) \tilde{s}(r')x$. Third, it is  a generalized bialgebra (and generalized weak bialgebra) in the sense that we have the axioms 
 $ \cop(x) (\tilde{t}(r) \o 1) = $ $\cop(x)(1 \o \tilde{s}(r))$,
$\cop(xy) = \cop(x)\cop(y)$ (which makes sense thanks to the previous axiom), $\cop(1) = 1 \o 1$,
and  $\eps(1_{S'}) = 1_{R'}$ and
 $\eps(xy) = \eps(x\tilde{s}(\eps(y)) ) = \eps(x\tilde{t}(\eps(y)) )$  for all $x, y \in T', r,r' \in R'$. 
A right bialgebroid is defined like a left bialgebroid with three of the axioms transposed \cite{KS}. 

 In \cite[4.1]{KS} we establish that $S := \End {}_BA_B$ is a left bialgebroid over $R$ with the $R$-$R$-bimodule
structure given by $$ r \cdot \alpha \cdot r' := \lambda(r) \rho(r') \alpha = r \alpha(-)r'$$
 for $r,r'$ in the centralizer $R$.  
The comultiplication $\cop_S: S \to S\o_R S$ is given by either of two formulas: 
\begin{equation} 
\label{eq:comult}
\cop_S(\alpha) := \sum_i \alpha(-t_i^1)t_i^2 \o \beta_i = \sum_j \gamma_j \o u_j^1 \alpha(u_j^2 -). 
\end{equation}
Since $S \o_R S \cong \Hom_{B-B}(A \o_B A, A)$ via $\alpha \o \beta \mapsto$ $ (a \o a' \mapsto \alpha(a)\beta(a'))$
\cite[3.11]{KS}, our formulas for the coproduct simplify greatly via eqs.~(\ref{eq: d2 quasibase})
to $\cop_S(\alpha)(a \o a') = \alpha(aa')$ (a generalized Lu bialgebroid).  The counit belonging to this coproduct
is $\eps_S: S \to R$ given by $\eps_S(\alpha) = \alpha(1_A)$.  It is now easy to see
that both maps are $R$-bimodule morphisms, that $\cop_S$ is coassociative and that $(\eps_S \o \id)\cop_S =$ $ \id_S = $ $(\id_S \o \eps_S) \cop_S$.

The ring $T$ defined above for any ring extension is a right bialgebroid over the centralizer 
$R$:  $S$ and $T$ are dual bialgebroids w.r.t.\ either of the nondegenerate pairings $\bra \alpha \| t \ket := \alpha(t^1)t^2$ or $[ \alpha \| t ] := t^1 \alpha(t^2) $, both
with values in $R$ \cite[5.3]{KS}.  The $R$-coring structure underlying the right $R$-bialgebroid $T$ is given by $r \cdot t \cdot r' = $ $ rt^1 \o t^2 r'$,
\begin{equation}
\label{eq: tee}
\cop_T(t) := \sum_j (t^1 \o_B \gamma_j(t^2)) \o_R u_j = \sum_i t_i \o_R (\beta_i(t^1) \o_B t^2)
\end{equation}
and counit $\eps_T(t) = t^1 t^2$, the multiplication mapping restricted to $T$.
 
In \cite[4.1]{KS} we observed that $S$ acts on $A$ 
(via evaluation) as a left $S$-module algebra (or algebroid): if $A_B$ is a balanced module,
then the invariant subring $A^S = B$.  In this paper, we will be more concerned with the dual
concept, comodule algebra (defined below).  As an example of this duality, and a guide
to what we are about to do, we dualize, as we would (but more carefully)
 for  Hopf algebra actions, 
the left action just mentioned $\lact: S \o_R A \to A$,  $\alpha \lact a := \alpha(a)$ for $\alpha \in S, a \in A$, to a right coaction
$\varrho_T: A \to A \o_R T$ given by $\rho(a) = a\0 \o a\1$ where
$\alpha \lact a = a\0 [\alpha | a\1]$. This comes out as $\varrho_T(a) = 
\sum_j \gamma_j(a) \o u_j$, since $\alpha(a) =$ $ \sum_j \gamma_j(a) [\alpha \| u_j]$  (obtained by applying
$\id \o \alpha$ to eq.~(\ref{eq: d2 quasibase})). The resulting right $T$-comodule algebra structure on $A$
is studied in \cite{Karl}, where it is shown that $A_B$ balanced results in a Galois extension $A \| B$ 
in the usual Galois coaction picture.\footnote{We expect the duality  left module algebra $\ \leftrightarrow \ $
right comodule algebra \newline
for f.g.\ projective bialgebroids to lead to duality for the notions
of Galois action \cite{Sz} and Galois coaction \cite{Karl,Bo}
for finite projective extensions 
(in \cite{BaS}?).
Indeed for the case of a depth two extension,
$A \| B$ is a Galois extension from both points of view \cite{Sz,Karl}.}

There is also an action of $T$ on $\ME$ studied 
in \cite[5.2]{KS}: the $R$-bialgebroid $T$ acts from the right on $\ME$ by
$f \ract t := $  $ f(-t^1)t^2$ for $f \in \ME, t \in T$.
This action makes $\ME$ a right $T$-module algebra
with invariants $\rho(A)$ (where $\rho(a)(x) = xa$ for $x,a \in A$). Thinking in terms of Hopf algebra
duality, we then expect to see a left coaction $\varrho: \ME \to S \o_R \ME$ with Sweedler
notation $\varrho(f) = f\-1 \o f\0$ 
satisfying $f \ract t = [f\-1 | t] f\0$. This comes in terms of a right D2 quasibase as 
\begin{equation}
\label{eq: row}
\varrho(f) = \sum_j \gamma_j \o (f \ract u_j)
\end{equation}
since $f \ract t = $ $\sum_j t^1 \gamma_j(t^2)u_j^1 f(u_j^2-) = $
$ \sum_j [\gamma_j \| t ] (f \ract u_j)$.  Since
$\ME$ is a variant of a smash product of $A$ with $S$
(cf.\ \cite[section 3]{KS}), we would want to 
show that $\ME$ is a Galois extension of a copy of $A^{\rm op}$
somewhat in analogy with cleft Hopf algebra coaction
theory although there is no antipode in our set-up:  see eq.~(\ref{eq: total integral}) for why
we may think of the natural inclusion $S \into \ME$ as
a total integral which cleaves the $S$-extension $\ME \| \rho(A)$.    
We next turn to several definitions,
lemmas and a theorem below in which we prove that $\varrho$ is a Galois coaction
for the left $S$-extension $\ME$ over $\rho(A)$.

\begin{definition}
Let $S'$ be a left $R'$-bialgebroid $(S', \tilde{s}, \tilde{t},$ $ \cop, $ $ \eps)$.   
A (left) $S'$-comodule algebra $C$ is a ring $C$ with ring homomorphism $R' \to C$ 
 together with a coaction $\delta: C \to S' \o_{R'} C$, where  values $\delta(c)$ are denoted by the Sweedler
notation $c\-1 \o c\0$, 
such that $C$ is a left $S'$-comodule over the $R'$-coring $S'$ \cite[18.1]{BW}, 
\begin{equation}
 \delta(1_C) = 1_{S'} \o 1_C, 
\end{equation} 
\begin{equation}
\label{eq: homog}
c\-1 \tilde{t}(r) \o c\0 = c\-1 \o c\0 \cdot r
\end{equation}
for all $r \in R$,
and 
\begin{equation}
\label{eq: cop homo}
\delta(cc') = \delta(c) \delta(c')
\end{equation}
 for all $c,c' \in C$.   The subring of coinvariants
is \begin{equation}
{}^{\rm co \, S'}C := \{ c \in C \| \delta(c) = 1_{S'} \o c \}.
\end{equation} 
$C$ is said to be a left $S'$-extension of ${}^{\rm co \, S'}C $. 
\end{definition} 

Like in the definition of left bialgebroid, the axiom~\ref{eq: cop homo} makes sense
because of the axiom~\ref{eq: homog}. 
The ring homomorphism $R' \to C$ induces a natural $R'$-bimodule
on $C$ which we refer to implicitly.

\begin{definition}
Let $S'$ be a right finite projective left $R'$-bialgebroid.  
An $S'$-comodule algebra $C$ is a left $S'$-Galois extension of its coinvariants
$D$ if the (Galois) mapping 
$\beta: C \o_D C \to S' \o_{R'} C$ defined
by $\beta(c \o c') = c\-1 \o c\0 c' $ is bijective.\footnote{These definitions are the left-handed versions of right comodule algebra
and right Galois extension for right bialgebroids in \cite[2.4,2.5]{LK} (cf.\ \cite[31.23]{BW}).  The definition above
for Galois extension is equivalent to $S' \o_{R'} C$ being a Galois $C$-coring \cite{BW},
an approach taken in \cite{Karl}, but not pursued in the present paper.}  
\end{definition} 
 
We need several lemmas for the nontrivial task of proving $\ME \| \rho(A)$ a left Galois extension.
The next lemma will be used among other things to show that $\varrho$ is coassociative. 

\begin{lemma}
\label{lemma: a}
Let $A \| B$ be a D2 extension. Then we have the isomorphisms
\begin{equation}
\label{eq: short iso}
S \o_R \ME \cong \Hom ({}_BA \o_B A, {}_B A) 
\end{equation}
via $\alpha \o f \longmapsto (a \o a' \mapsto \alpha(a)f(a'))$, and
\begin{equation}
\label{eq: long iso}
S \o_R S \o_R \ME \cong \Hom ({}_BA \o_B A \o_B A, {}_B A)
\end{equation}
via $\alpha \o \beta \o f \longmapsto (a \o a' \o a'' \mapsto
\alpha(a) \beta(a') f(a''))$.
\end{lemma}
\begin{proof}
The inverse in (\ref{eq: short iso}) is given by $F \mapsto $
$ \sum_j \gamma_j \o u^1_j F(u^2_j \o -) $ by eq.~\ref{eq: d2 quasibase}.

The inverse in (\ref{eq: long iso}) is given by $$F \mapsto 
\sum_{j,k,i} \gamma_j \o \gamma_k \o u^1_k u^1_j F(u^2_j \gamma_i(u^2_k)u^1_i \o u^2_i \o -) 
$$ since
$$ \sum_{j,k,i} \gamma_j(a)\gamma_k(a')u^1_k u^1_j F(u^2_j \gamma_i(u^2_k)u^1_i \o u^2_i \o a'') 
 = \sum_{j,k} \gamma_j(a)u^1_jF(u^2_j \gamma_i(a') u^1_i \o u^2_i \o a'')$$
$ = F(a \o a' \o a'') $
and
$$ \sum_{j,k,i} \gamma_j \o \gamma_k \o_R u^1_k u_j^1 \alpha(u^2_j \gamma_i(u^2_k)u^1_i) \beta(u^2_i) f(-) =
\sum_{j,i} \gamma_j \o u^1_j \alpha(u^2_j \gamma_i(-)u^1_i) \beta(u^2_i) \o f(-) $$
$ = \sum_j \gamma_j(-) u^1_j \alpha(u^2_j) \o \beta \o f = \alpha \o \beta \o f$
for $f \in \ME, \alpha, \beta \in S$. 
\end{proof}

The existence alone of an isomorphism in the next lemma may
be seen by letting $M$ be free of rank one. 

\begin{lemma}
Given any rings $A$ and $B$, with modules $M_A$, ${}_BU$ and bimodule ${}_BN_A$ with $M_A$ f.g.\ projective,
then 
\begin{equation}
\label{eq: mull}
M \o_A \Hom({}_B N, {}_B U) \stackrel{\cong}{\longrightarrow} \Hom ({}_B\Hom (M_A, N_A), {}_B U) 
\end{equation}
via the mapping $m \o \phi \longmapsto (\nu \mapsto \phi (\nu(m)))$.
\end{lemma}
\begin{proof}
Let $m_i \in M$, $g_i \in \Hom (M_A, A_A)$ be dual bases for $M_A$.
For each $n \in N$ let $ng_i(-)$ denote the obvious mapping
in $\Hom (M_A, N_A)$.  
Then the inverse mapping is given by 
\begin{equation}
F \longmapsto \sum_i m_i \o (n \mapsto F(ng_i(-))).
\end{equation}
  Note that both maps are well-defined module homomorphisms,
and inverse to one another since 
$$ \sum_i m_i \o \phi(-g_i(m)) = \sum_i m_i \o_A (g_i(m) \phi)(-) =
\sum_i m_i g_i(m) \o \phi = m \o \phi$$
and $ \nu \mapsto \sum_i F(\nu(m_i) g_i ) = F(\nu)$ for
$\nu \in $ $\Hom (M_A, N_A)$. 
\end{proof}
The lemma above is relevant in our situation since the depth two condition
implies that a number of constructions such as the tensor-square and 
endomorphism rings are finite projective. For example, $\ME_A$ is f.g.\ 
projective \cite[3.13]{KS}, which we may also see directly from eq.~(\ref{eq: d2 quasibase})
by applying $\id_A \o_B f$ for $f \in \ME$, viewing $\gamma_j \in S \subseteq \ME$
and an obvious mapping  of $A \o_B A$ into $\Hom (\ME_A, A_A)$ which appears in the next
lemma.

\begin{lemma}
\label{lemma: b}
If $A \| B$ is D2, then
$${}_B A \o_B A \stackrel{\cong}{\longrightarrow} {}_B \Hom(\ME_A, A_A) $$
via $\Psi(a \o a')(f) = af(a') $.
\end{lemma}
\begin{proof}
 Let $F \in (\ME_A)^*$ (the right $A$-dual of $\ME$).  Define an inverse 
$\Psi^{-1}(F) = \sum_j  F(\gamma_j)u_j$
 where $\gamma_j \in S, u_j \in T$ is a right D2 quasibase.  Then $\Psi^{-1} \Psi = \id_{A \o_B A}$ by eq.~(\ref{eq: d2 quasibase}).
Also $\Psi \Psi^{-1} = \id_{\ME^*}$ since $$\Psi \Psi^{-1}(F)(f) = \sum_j F(\gamma_j)u^1_j f(u^2_j) = F(\sum_j \gamma_j(-) u^1_j f(u^2_j)) = F(f). \qed $$ 
\renewcommand{\qed}{}\end{proof}

Recall that $\rho(A)$ denotes the set of right multiplication operators by elements of $A$.  

\begin{theorem}
\label{th-endogalext}
Let $A \| B$ be a depth two extension.  Then $\ME$ is a left $S$-comodule algebra with the coaction~(\ref{eq: row})
and a Galois extension of its coinvariants $\rho(A)$.  
\end{theorem}

\begin{proof}
Recall that the coaction $\varrho$ is given on $f \in \ME = \End {}_BA$  by $$ f\-1 \o f\0 = \sum_j \gamma_j \o u_j^1 f(u_j^2 -)$$
where $\gamma_j \in S, u_j \in T$ is a right D2 quasibase. First, the ring homomorphism $R \to \ME$ is 
given by $\lambda$, so 
for $1_{\ME} = \id_A = 1_S$, we have 
$$ \varrho(\id_A) = \sum_j \gamma_j \o_R \lambda(u^1_j u^2_j) = \sum_j \gamma_j(-)u_j^1 u_j^2 \o \id_A = 1_S \o 1_{\ME}.
$$

Secondly, we check that $\ME$ forms a left $S$-comodule w.r.t.\ the $R$-coring $S$ and the coaction $\varrho$. 
The coaction is coassociative, $(\cop_S \o \id_{\ME})\varrho =$ $ (\id_S \o \varrho)\varrho$,
for we use lemma~\ref{lemma: a} (as an identification and suppressing the isomorphism) and eqs.~(\ref{eq:comult})
and~(\ref{eq: tee mult})  to check values of each side of this equation, evaluated on $A \o_B A \o_B A$:
\begin{eqnarray*}
 (\sum_j \cop_S(\gamma_j) \o (f \ract u_j))(a \o a' \o a'') &=& \sum_{k,j} \gamma_k(a)u_k^1\gamma_j(u^2_k a')u_j^1 f(u_j^2 a'')  \\
& = & \sum_j \gamma_j(aa') u_j^1 f(u^2_ja'') = f(aa'a'') \\
& = & \sum_{i,j} \gamma_i(a) \gamma_j(a')u^1_j u^1_i f(u^2_i u^2_j a'') \\
& = & (\sum_j \gamma_j \o \varrho (f \ract u_j))(a \o a' \o a''). 
\end{eqnarray*}
We must also check that $\varrho: \ME \to S \o_R \ME$ is a left $R$-module morphism: given $r \in R, f\in \ME$, we
 use lemma~\ref{lemma: a} again to note that for $a, a' \in A$ 
$$ \varrho(\lambda(r)f) (a \o a') = \sum_j \gamma_j(a)u^1_j rf(u^2_j a') = rf(aa') = (r\cdot f\-1 \o f\0)(a \o a'), $$
by an application of eq.~(\ref{eq: d2 quasibase}) (inserting an $r \in C_A(B)$).  
(Note with $r = 1$ that we obtain  
\begin{equation}
\label{eq: total integral}
\varrho(\alpha) = \cop_S(\alpha) \ \ \ \ \ \forall \alpha \in S,
\end{equation}
which should be compared to the total integral and cleft extension approach in \cite[6.1]{KN}.)
Finally, $\ME$ is counital, whence a left $S$-comodule, since for $f \in \ME$, 
$$ (\eps_S \o \id_S)\varrho(f) = \sum_j \gamma_j(1)u^1_jf(u^2_j -) = f. $$

Next we must check that $\Im \varrho$ lies in a submodule of $S \o_R \ME$ where tensor product multiplication
makes sense: again using lemma~\ref{lemma: a} and for $a,a' \in A$,
$$ (f\-1 \tilde{t}(r) \o f\0)(a \o a') = \sum_j \gamma_j(ar)u_j^1 f(u_j^2 a') = f(ara') $$
$$ = \sum_j \gamma_j(a) u_j^1 f(u_j^2 ra') = (f\-1 \o f\0 \lambda(r))(a \o a'). $$
Then multiplicativity of the coaction follows from the measuring axiom satisfied by the right action of
$T$ on $\ME$ \cite[5.2]{KS} and eq.~(\ref{eq: tee}) ($f,g \in \ME$):  
$$ \varrho(fg)(a \o a') = \sum_j (\gamma_j \o (f\ract {u_j}\1) \circ (g \ract {u_j}\2))(a \o a') = $$
$$ \sum_{j,k} \gamma_j(a) u^1_j f(\gamma_k(u^2_j) u^1_k g(u^2_k a')) = f(g(aa')) = \sum_{j,i} \gamma_i(\gamma_j(a))
u^1_i f(u^2_i u^1_j g(u^2_j a')) = $$
$$ = \sum_{i,j} (\gamma_i \circ \gamma_j \o (f \ract u_i ) \circ (g \ract u_j))(a \o a') = \varrho(f)\varrho(g)(a \o a'). $$

Next we determine the coinvariants ${}^{\rm co \, S} \ME$. Given $a \in A$,
we note that $\rho(a) \in {}^{\rm co \, S} \ME$ since
$$ \varrho(\rho(a)) = \sum_j \gamma_j \o u^1_j u^2_j(-a) = 1_S \o \rho(a).$$
Conversely, suppose $\sum_j \gamma_j \o (f \ract u_j) = 1_S \o f$
in $S \o_R \ME \cong \Hom ({}_BA \o_B A, {}_BA)$,
then for $a, a' \in A$,
$$ \sum_j \gamma_j(a)u^1_j f(u^2_j a') = f(aa') = af(a'). $$
It follows that $f(a) = af(1_A)$ for all $a \in A$, so
$f = \rho(f(1)) \in \rho(A)$.  Hence, ${}^{\rm co \, S} \ME = \rho(A)$. 

Finally, the Galois mapping $$\beta: \ME \o_{\rho(A)} \ME \to
S \o_R \ME, \ \ \ \ \ \beta(f \o g) = f\-1 \o f\0 g $$
under the identification $S \o_R \ME \cong$ 
$\Hom ({}_B A \o_B A, {}_BA)$ in lemma~\ref{lemma: a}
is given by an application of eq.~(\ref{eq: d2 quasibase}):  ($a,a' \in A, f,g \in \ME$)
 
\begin{equation}
\label{eq: galois map}
\beta(f \o g)(a \o a') = \sum_j \gamma_j (a)u^1_j f(u^2_j g(a')) = f(ag(a')).
\end{equation}

We show $\beta$ to be a  composite of several isomorphisms using the lemmas
(commutative diagram in section~1).  First note
that 
$\rho(A) \cong  A^{\mbox{op}}$ 
and 
${}_{\rho(A)} 
\ME_{\rho(A)}$
 given by
 $\rho(a')\circ f \circ \rho(a) (a'') = f(a'' a) a'$ 
is equivalent to ${}_A \ME_A$ given by $a \cdot f \cdot a' (a'') = f(a'' a)a'$.  
This is the usual $A$-$A$-bimodule structure on the left endomorphism ring $\ME$
considered in \cite[3.13]{KS}, where  $\ME_A$ is shown to be f.g.\ projective.
Consider then the composition of isomorphisms,
\begin{eqnarray*}
 \ME \o_{\rho(A)} \ME  & \stackrel{\cong}{\longrightarrow} & \ME \o_A \Hom ({}_BA, {}_BA) \\
& \stackrel{\cong}{\longrightarrow} & \Hom ({}_B \Hom(\ME_A, A_A), {}_BA) \\
& \stackrel{\cong}{\longrightarrow} & \Hom ({}_B A \o_B A, {}_B A) 
\end{eqnarray*}
given by $$ f \o g \mapsto g \o f \mapsto (\nu \mapsto f(\nu(g))) \mapsto ( a\o a' \mapsto  f( a g(a')). $$
This is $\beta$ as given in eq.~(\ref{eq: galois map}), whence 
$\beta$ is an isomorphism and the extension $\ME \| \rho(A)$ is  Galois.  
\end{proof}

\section{A Hopf algebroid for all seasons}

One of the few regrets of generalizing Lu's bialgebroid $A^e$ over an
algebra $A$ to the right bialgebroid $T$ of a depth two extension
$A \| B$ is that the antipode is lost, for the flip or twist anti-automorphism on $A \o_k A^{\rm op}$ 
does not extend to a self-mapping of $(A \o_B A)^B$. If we require $B$ to be separable with symmetric separability
element however, there is a projection of $A \o_B A \to T$ which we may apply to define
a twist of $T$.  However, 
Lu's definition of Hopf algebroid \cite{Lu} makes it necessary to find an appropriate section $T \o_R T \to T \o_k T$ of the canonical
map in the other direction, although  the centralizer $R$
is not \textit{a priori} separable.  In this section, we carry out 
this plan using instead the alternative definition
of Hopf algebroid in \cite[B\"{o}hm and Szlach\`{a}nyi]{BS}.    
 
Let $k$ be a commutative ring and $B$ a Kanzaki separable $k$-algebra
\cite[strongly separable algebra]{KSt}.  This means that there is a separability element
$e = e^1 \o e^2$ $ \in B \o_k B$ which is symmetric, so
that $e^1 e^2 = 1 $ $ = e^2 e^1$ as well as $be = eb$
and $e^1 b \o e^2 = $ $ e^1 \o be^2$ for all $b \in B$. (Typically
for quantum algebra, we use both
these equalities repeatedly below together with $bt = tb$ for $t \in T$ and well-definedness of various mappings for commuting elements in the equations in the proof of the
theorem.)
For example, all separable algebras over a field of characteristic zero
are Kanzaki separable.  Over a field of characteristic $p$ matrix
algebras of order divisible by $p$ are separable although not
Kanzaki separable.  Fix the notation above for the next theorem.

\begin{theorem}
Let $A \| B$ be a depth two extension of $k$-algebras where
$B$ is Kanzaki separable. Then the left bialgebroid
$T^{\rm op}_{\rm cop}$ is a Hopf algebroid.
\end{theorem}
\begin{proof}
The standard right bialgebroid $T = (A \o_B A)^B$ with structure $(T, R, s_R, t_R,$ $ \cop, \eps)$ becomes a left
bialgebroid via the opposite multiplication and coopposite comultiplication as follows \cite[2.1]{KS}:  
$$T^{\rm op}_{\rm cop} := ( T^{\rm op}, R^{\rm op}, s_L = s_R, t_L = t_R, \cop^{\rm op}, \eps). $$   
The multiplication on $T^{\rm op}$ is 
\begin{equation}
tt' = t^1 {t'}^1 \o_B {t'}^2 t^2 
\end{equation}
while the target and source maps are $s_L: R^{\rm op} \to T^{\rm op}$, $t_L: R \to T^{\rm op}$ are then
$s_L(r) = 1_A \o r$, $t_L(r) = r \o 1_A$ for $r \in R$.  The $R^{\rm op} $ - $R^{\rm op}$-bimodule structure
on $T$ is then given by
\begin{equation}
r \cdot t \cdot r' = s_L(r) t_L(r') t = (r' \o r)t = r' t^1 \o t^2 r 
\end{equation}
in other words ${}_{R^{\rm op}}T_{R^{\rm op}}$ is the standard bimodule ${}_RT_R$ with endpoint multiplication after passing
to modules over the opposite algebra of $R$.  Tensors over $R^{\rm op}$ are the same as tensors over $R$ after a flip;
e.g., $$
T \o_{R^{\rm op}}T \stackrel{\cong}{\longrightarrow} T \o_R T \stackrel{\cong}{\longrightarrow} (A \o_B A \o_B A)^B 
$$
via the mapping 
\begin{equation}
\label{eq: iq}
t \o t' \stackrel{\cong}{\longmapsto} {t'}^1 \o {t'}^2 t^1 \o t^2,
\end{equation}
which is an $R^{\rm op}$-$R^{\rm op}$-isomorphism \cite[5.1]{KS}. 

The comultiplication $\cop^{\rm op}: T \to T \o_{R^{\rm op}} T$ is given for $t \in T$ by
the coopposite of eq.~(\ref{eq: tee}), 
\begin{equation}
\label{eq: tee cop}
 \cop^{\rm op}(t)  = \sum_j u_j \o (t^1 \o_B \gamma_j(t^2))
\end{equation}
which in $(A \o_B A \o_B A)^B$ is the value $t^1 \o 1_A \o t^2$ after applying eq.~(\ref{eq: iq})
and the right D2 quasibases eq.~(\ref{eq: d2 quasibase}). We will denote below $\cop^{\rm op}(t) = $ $t\1 \o t\2$. 

The antipode $\tau: T \to T$ is a flip  composed with a projection from $A \o_k A$:
\begin{equation}
\label{eq: antipode}
\tau(t) := e^1 t^2 \o_B t^1 e^2.
\end{equation}
 We next note that $\tau$ is an anti-automorphism of order two on $T^{\rm op}$. Let $e = f$ in $B \o_k B$ so that  
$$ \tau(t') \tau(t) = e^1 {t'}^2 f^1 t^2 \o_B t^1 f^2 {t'}^1 e^2 = e^1 {t'}^2 t^2 \o_B t^1 f^2 f^1 {t'}^1 e^2 = $$
$ = \tau(tt')$ since ${t'}^1 e^2 \o_k e^1 {t'}^2 \in (A \o_k A)^B$ and $f^2 f^1 = 1_B$. In addition,
$$ \tau^2(t) = \tau(e^1 t^2 \o t^1 e^2) = f^1 t^1 e^2 \o_B e^1 t^2 f^2 = t, $$
since $e^2 e^1 = 1$.

Next, we show that $\tau$ satisfies the three axioms of \cite[Def.\ 4.1]{BS} given below in
eqs.~(\ref{eq: (4.1)})-(\ref{eq: (4.3)}).  Note that $\tau^{-1} = \tau$.  First,
\begin{equation}
\label{eq: (4.1)}
\tau \circ t_L = s_L
\end{equation}
since $\tau(r \o 1) = e^1 \o_B r e^2 = 1 \o r$ for $r$ in the
centralizer $C_A(B)$.  

Second, we have the equality in $T\o_{R^{\rm op}} T$, 
\begin{equation}
\label{eq: (4.2)}
\tau^{-1}(t\2)\1 \o \tau^{-1}(t\2)\2 t\1 = \tau^{-1}(t) \o 1_T
\end{equation}
which follows from eq.~(\ref{eq: tee cop}) and (working from left to right): 
$$ \sum_{j,k} u_k \o (e^1 \gamma_j(t^2)  \o_B  \gamma_k(t^1 e^2))u_j =
\sum_{j,k} u_k \o (e^1 \gamma_j(t^2)u_j^1 \o_B  u_j^2\gamma_k(t^1 e^2)) $$
the last expression mapping via the isomorphism in eq.~(\ref{eq: iq}) into
$$e^1 \o_B t^2 \gamma_k(t^1) u_k^1 \o_B u^2_k e^2 = 1_A \o e^1 t^2 \o t^1 e^2 $$
in $(A \o_B A \o_B A)^B$, which is the image of $\tau^{-1}(t) \o 1_T$ under the
same isomorphism.  

Finally, we have the equality in $T\o_{R^{\rm op}} T$,
\begin{equation}
\label{eq: (4.3)}
\tau(t\1)\1 t\2 \o \tau(t\1)\2 = 1_T \o \tau(t)
\end{equation}
for all $t \in T$, which follows similarly from
$$\sum_{j,k} u_k(t^1 \o_B \gamma_j(t^2)) \o e^1 u^2_j \o_B \gamma_k(u^1_j)e^2) = 
\sum_{j,k} u_k^1 t^1 \o_B \gamma_j(t^2)u_k^2) \o e^1 u^2_j \o_B \gamma_k(u^1_j)e^2)$$
which maps isomorphically to
$$\sum_{j,k} e^1 u^2_j \o_B \gamma_k(u^1_j)u^1_k t^1 \o_B \gamma_j(t^2) u^2_k e^2 = \sum_j
e^1 u^2_j \o t^1 \o \gamma_j(t^2) u^1_j e^2 = $$
$$ e^1 t^2 \o t^1 e^2 \o 1_A \stackrel{\cong}{\longleftarrow} 1_T \o \tau(t) $$
since $\sum_j \gamma_j(a)u^1_j e^2 \o_k e^1 u^2_j = e^2 \o_k e^1a$ for $a \in A$
follows from eq.~(\ref{eq: d2 quasibase}). 
\end{proof}
  
Suppose $k$ is a field, then $B$ is finite dimensional as it is a separable algebra. If moreover 
 $A \| B$ is a proper ring extension, then there is bimodule projection $A \to B$ by separability,
whence $A_B$ is finitely generated (and projective \cite[2.2]{Karl}) and so $A$ is finite dimensional as well.  
The theorem is thus viewed as a natural generalization of Lu's twist Hopf algebroid to certain finite-dimensional algebra pairs.     

\section{Galois extended Hopf subalgebras are normal}

There is a question of whether depth two Hopf subalgebras are normal \cite[3.4]{KK}. 
In this section we answer this question in an almost unavoidable
special case, namely, when the Hopf subalgebra forms a Hopf-Galois
extension with respect to the action of a third Hopf algebra.  
Since a depth two extension with one extra condition is a Galois extension for
actions of bialgebroids or weak bialgebras \cite{LK}, 
the situation of
ordinary Hopf-Galois extension would seem to be a critical  step.  

Let $k$ be a field. All Hopf algebras in this section are finite dimensional
algebras over $k$.  
Recall that a Hopf subalgebra $K \subseteq H$ is a Hopf algebra $K$
w.r.t.\  the algebra and coalgebra structure of $H$ (with counit denoted by $\eps$) as well as
stable under the antipode $\tau$ of $H$. Recall the Nichols-Zoeller
result that the natural modules $H_K$ and ${}_KH$ are free.  
$K$ is \textit{normal} in $H$
if $\tau(a\1)x a\2 \in K$ and $a\1 x \tau(a\2) \in K$
for all $x \in K, a \in H$.  Equivalently, if $K^+$ denotes
the kernel of the counit $\eps$, $K$ is a normal Hopf subalgebra
of $H$ if the left algebra ideal and coideal $HK^+$ is equal to the right ideal and coideal $K^+H$
\cite[3.4.4]{Mo}.
   
In considering another special case of D2 Hopf subalgebras,  we showed in \cite{Karl} that H-separable Hopf subalgebras
are normal using favorable properties for H-separable extensions of going down and going up for ideals.  However, we
noted that
such subalgebras are not proper if $H$ is semisimple, e.g.,
$H$ is a complex group algebra.
In \cite[3.1]{KK} we showed that depth two subgroups are normal subgroups
using character theory (for $k = \C$). We also noted the more general
converse that normal Hopf $k$-subalgebras are Hopf-Galois extensions
and therefore D2. Next we extend this to the 
characterization of normal Hopf subalgebras below,  one that we believe 
is not altogether unexpected but unnoted or not adequately exposed in the literature.

  \begin{theorem}
Let $K \subseteq H$ be a Hopf subalgebra.  Then $K$ is normal in
$H$ if and only if $H | K$ is a Hopf-Galois extension.  
\end{theorem}

\begin{proof} ($\Rightarrow$) This is more or less implicit
in \cite[3.4.4]{Mo}, where it is also shown \cite[chs.\ 7,8]{Mo}
 that $H$ is a crossed product by a counital $2$-cocycle of $K$ with 
the quotient Hopf algebra $\overline{H}$ (a cleft $\overline{H}$-extension or
Galois extension with normal basis).
Since $HK^+ = K^+H$ under normality of $K$, it becomes a Hopf ideal, so 
we form the Hopf algebra $\overline{H} := H/HK^+$, which coacts
naturally on $H$ via the comultiplication and quotient projection.
The coinvariants are precisely $K$ since $H_K$ is faithfully flat.
The Galois mapping $\beta: H \o_K H \rightarrow H \o \overline{H}$
given by $\beta(a \o a') = a {a'}\1 \o \overline{{a'}\2}$ is an isomorphism
with inverse given by $x \o \overline{y} \mapsto x \tau(y\1) \o y\2$. 

($\Leftarrow$)  Suppose $H$ is a $W$-Galois extension of $K$ where $W$ is a Hopf algebra
with right coaction $\rho: H \to H \o W$ on $H$.  We define a mapping $\Phi: H \to W$
by $\Phi(h) = \eps_H(h\0) h\1$, i.e., $\Phi = (\eps_H \o \id_W) \circ \rho$.  We note
that $\Phi$ is an algebra homomorphism since $\rho$ and $\eps_H$ are (and augmented since $\eps_W \circ \Phi = \eps_H$). 
Also, $\Phi: H \to W$ is a right $W$-comodule morphism since $H$ is a right $W$-comodule with
$\rho$ and $\cop_W$ obeying a coassociativity rule.  Next we note that $\Phi$ is an epi
since given $w \in W$, there is $\sum_i h_i \o {h'}_i \in H \o_K H$ such that
$1 \o w = \sum_i h_i {{h'}_i}\0 \o {{h'}_i}\1$.  Applying $\eps_H \o \id_W$ to this, we obtain 
$$ \Phi(\sum_i \eps_H(h_i) {h'}_i) = w . $$

We note that $\ker \Phi$ contains $K^+$ since $K = H^{\rm co \, W}$ $ = \{ h \in H \, |\, \rho(h) = h \o 1_W \}$
Consider the coalgebra and right quotient $H$-module $H/ K^+H := \overline{H}$
as well as the coalgebra and left quotient $H$-module $H/ HK^+ := \overline{\overline{H}}$.
In this case, $\Phi$ induces $\overline{\Phi}:
\overline{H} \to W$ and $\overline{\overline{\Phi}}: \overline{\overline{H}} \to W$.
(They are respectively right and left $H$-module morphisms w.r.t.\ the modules $W_{\Phi}$ and ${}_{\Phi}W$.)
By Schneider \cite[1.3]{HJS}, the Galois quotient mapping $\overline{\beta}: H \o_K H \to H \o \overline{H}$
given by $\overline{\beta}(x \o y) = xy\1 \o \overline{y\2}$ is bijective (since $K$ is a left coideal subalgebra of $H$).
But the Hopf subalgebra $K$ is also a right coideal subalgebra satisfying a right-handed version of Schneider's lemma
recorded in \cite[2.4]{FMS}: whence the Galois mapping $\overline{\overline{\beta}}: H \o_K H \to \overline{\overline{H}} \o H$ given by $\overline{\overline{\beta}}(x \o y) = \overline{x\1} \o x\2 y$ is bijective as well.

Observe now that $H_K$ is free of rank $n$, let's say,  so $\beta$ bijective implies that $\dim_k W = n$.
Similarly, $\overline{\beta}$ bijective implies $\dim_k \overline{H} = n$
and $\overline{\overline{\beta}}$ bijective implies $\dim_k \overline{\overline{H}} = n$.  
It follows that the  vector space epimorphisms $\overline{\Phi}: \overline{H} \to W$ and $\overline{\overline{\Phi}}: \overline{\overline{H}} \to W$ are isomorphisms. But $\overline{\Phi}$ factors through $\overline{H} \to H/HK^+H$ induced
by $K^+H \subseteq HK^+H$; similarly, $\overline{\overline{\Phi}}$ factors through $\overline{\overline{H}} \to
H/HK^+H$, so both these canonical mappings are monic.  
It follows that $HK^+ = HK^+H = K^+H$, whence $K$ is a normal Hopf subalgebra in $H$. 
  \end{proof}

In the proof of $\Leftarrow$ above, we can go further to conclude that $\overline{H}$ is a Hopf algebra
isomorphic to $W$ as augmented algebras.  However, the theory of deforming the comultiplication of
a Hopf algebra by a $2$-cocycle
\cite[2.3.4]{Maj} shows that there are pairs of Hopf algebras isomorphic as augmented algebras yet
non-isomorphic as Hopf algebras. Additionally, there are examples of Hopf-Galois extensions w.r.t.\ two
different Hopf algebras.  We therefore do not know \textit{a priori}
 if $\overline{H}$ and $W$ are isomorphic as Hopf algebras.

\section{Weak Hopf-Galois extensions are depth two}

In this section we study right Galois extensions of special bialgebroids - the weak Hopf-Galois extensions, cf.\ \cite{BW,CDG, Karl, LK}. By exploiting the antipode in weak Hopf algebras, we find an alternative Galois mapping
which characterizes weak Hopf-Galois extensions.  This leads
to several corollaries that finite weak Hopf-Galois extensions are right as well as left depth two
extensions, that they may be defined by only a surjective Galois map, and that a weak Hopf algebra over
its target separable subalgebra is an example of such.  We propose a number of problems for further study
in the young subject of weak Hopf-Galois extensions.

  Weak Hopf algebras are a special case of Hopf algebroids - those 
with separable base algebra \cite{EN, KS}: the separable
algebra has an index-one Frobenius system which one uses to
convert mappings to the base and tensors over the base
to linear functionals and tensors over a ground field. There is an example of one step in how to conversely
view a weak Hopf algebra $H$ as a Hopf algebroid over its left coideal subalgebra $H^L$ in the proof of corollary~\ref{cor - olary} below.  

Let $k$ be a field. A weak Hopf algebra
$H$ is  first a weak bialgebra, i.e., a $k$-algebra and $k$-coalgebra $(H, \cop, \eps)$
such that the comultiplication $\cop: H \to H \o_k H$ is linear and multiplicative,
$\cop(ab) = \cop(a) \cop(b)$, and the counit is linear
just as for bialgebras;
however, the change (or weakening of the axioms) is that $\cop$ and $\eps$ may not be unital, $\cop(1) \neq
1 \o 1$ and $\eps(1_H) \neq 1_k$, but must satisfy 
\begin{equation}
1\1 \o 1\2 \o 1\3 = (\cop(1) \o 1)(1 \o \cop(1)) =
(1 \o \cop(1))(\cop(1) \o 1)
\end{equation}
 and $\eps$ may not be multiplicative, $\eps(ab) \neq \eps(a) \eps(b)$ but must satisfy $(a,b,c \in H$) 
\begin{equation}
\eps(abc) = \eps(ab\1) \eps(b\2 c) = \eps(a b\2) \eps(b\1 c).
\end{equation}
There are several important projections that result from these axioms:
\begin{eqnarray}
\Pi^L(x) &:=& \eps(1\1 x) 1\2 \\
\Pi^R(x) &:=& 1\1 \eps(x 1\2) \\
\overline{\Pi}^L(x) & := & 1\1 \eps(1\2 x) \\
 \overline{\Pi}^R(x) & := & \eps(x 1\1) 1\2 \ \ \ \ \ (\forall \, x \in H)
\end{eqnarray}
We denote $H^L := \Im \Pi^L = \Im \overline{\Pi}^R$
and $H^R := \Im \Pi^R = \overline{\Pi}^L$. (These subalgebras
are separable $k$-algebras in the presence of an antipode.)

In addition to being a weak bialgebra, a weak Hopf algebra
has an antipode $S: H \to H$ satisfying the axioms
\begin{eqnarray}
S(x\1)x\2 & = & \Pi^R(x) \label{eq: bns2} \\
x\1 S(x\2) & = & \Pi^L(x) \label{eq: bns1} \\
S(x\1) x\2 S(x\3) & = & S(x) \label{eq: bnsa} \ \ \ \  (\forall \, x \in H)
\end{eqnarray}
The antipode turns out to be bijective for finite dimensional
weak Hopf algebras (which we will assume for the rest of this
section), an anti-isomorphism of algebras with inverse
denoted by $\overline{S}$.  

The reader will note from the axioms above that a Hopf algebra is automatically a weak Hopf algebra.  For a weak Hopf algebra
that is not a Hopf algebra, consider 
a typical groupoid algebra such as $H = M_n(k)$, the $n \times n$-matrices
over $k$ (the groupoid here being a category with $n$ objects
where each $\Hom$-group has a single invertible arrow).
Let $e_{ij}$ denote the $(i,j)$-matrix unit.  
For example, $M_n(k)$ is a weak Hopf algebra
with the counit given by $\eps(e_{ij}) = 1$, comultiplication by $\cop(e_{ij})
= $ $e_{ij} \o e_{ij}$ and antipode given by $S(e_{ij}) = e_{ji}$
for each $i,j = 1,\ldots,n$
(extending the Hopf algebra structure of group algebras).  In this
case, $H^L = H^R$ and is equal to the diagonal matrices.  
The corresponding projections are given by
$\Pi^L(e_{ij}) = e_{ii}$ $= \overline{\Pi}^L(e_{ij})$ and $\Pi^R(e_{ij})= e_{jj}$
$= \overline{\Pi}^R(e_{ij})$. Note that $\eps(1_H) = n 1_k$ which is zero
if the characteristic of $k$ divides $n$.  

There are a number of equations in the subject that we
will need later (cf.\ \cite[2.8, 2.9, 2.24]{BNS}):
\begin{eqnarray}
\Pi^L & = & S \circ \overline{\Pi}^L \\
\Pi^R & = & S \circ \overline{\Pi}^R \label{eq: over2} \\
\overline{S}(a\2)a\1 & = & \overline{\Pi}^R(a) \label{eq: over1} \\
a\2 \overline{S}(a\1) &=& \overline{\Pi}^L(a) \\
a\1 \o \Pi^L(a\2) & = & 1\1 a \o 1\2  \\
\Pi^R(a\1) \o a\2 & = & 1\1 \o a 1\2  \label{eq: pi-are} \\
\Pi^R(a)b & = & b\1 \eps(ab\2)  \\
a\Pi^L(b) & = & \eps(a\1 b) a\2  \label{eq: pi-ell} \ \ \ \ \ (\forall \, a, b \in H) 
\end{eqnarray}
where e.g.\ eq.~(\ref{eq: over1}) follows from applying the inverse-antipode
to eqs.~(\ref{eq: over2}) and~(\ref{eq: bns2}). 

We recall the definition of a right $H$-comodule algebra $A$,
its subalgebra of coinvariants, and Galois coaction 
for $H$ a weak bialgebra (e.g.\ in \cite{CDG}):

\begin{definition}
Let $H$ be a weak bialgebra with $A, H$ both $k$-algebras.
$A$ is a right $H$-comodule algebra if there is a right $H$-comodule
structure $\rho: A \to A \o_k H$ such that $\rho(ab) = \rho(a) \rho(b)$ for each $a,b \in A$ and any of the equivalent
conditions for $\rho(a) := a\0 \o a\1$ are satisfied:
\begin{eqnarray}
1\0 \o 1\1 & \in & A \o H^L \label{eq: cdg2} \\\
a\0 \o \Pi^L(a\1) & = & 1\0 a \o 1\1 \label{eq: cdg1} \\
a\0 \o \overline{\Pi}^R(a\1) & = & a 1\0 \o 1\1 \ \ \ \ (\forall \, a\in H) \label{eq: cdg3} \\
1\0 \o 1\1 \o 1\2 & = & (\rho(1_A) \o 1_H)(1_A \o \cop(1_H)) 
\end{eqnarray}
The coinvariants are defined by 
$$ B := \{ b \in A \| b\0 \o b\1 = 1\0 b \o 1\1  = b 1\0 \o 1\1 \}, $$
the second equation following from equations directly above.
We say $A$ is a weak $H$-Galois extension of $B$ if
the mapping $\beta: A \o_B A \to A \o_k H$ given
by $\beta(a \o a') = a {a'}\0 \o {a'}\1$ is bijective onto
$$ \overline{A \o H} = (A \o H )\rho(1) = \{ a 1\0 \o h 1\1 \|
a \in A, h \in H \}. $$  
\end{definition}

For finite dimensional weak Hopf algebras and their actions,
we only need require $\beta$ be surjective in the definition of weak Hopf-Galois extension,
as $\beta$ is automatically injective by \cite{Bo, BTW} or corollary~\ref{cor-surj} below.  
Note that $\Im \rho \subseteq \overline{A \o H}$,
an $A$-$B$-sub-bimodule
and that $\beta$ is an $A$-$B$-bimodule morphism
w.r.t.\ the structure $a' \cdot (a \o h) \cdot b = a'ab \o h$
on $A \o H$.  
These definitions correspond to the case of a separable base algebra in the definitions of right 
comodule algebras, Galois coring and Galois coactions for bialgebroids
given in \cite{Karl, LK}.
 
We now establish the Hopf algebra analogue of an alternate 
Galois mapping characterizing Galois extension. This would
correspond to working with a left-handed version of
the Galois coring considered in \cite{CDG}. 

\begin{prop}
Suppose $H$ is a weak Hopf algebra and $A$ a right $H$-module
algebra with notation introduced above.  
Let $\beta': A \o_B A \to A \o H$ be defined by 
\begin{equation}
\beta'(a \o a') = a\0 a' \o a\1
\end{equation}
and $\eta: A \o H \to A \o H$ be the map defined by
\begin{equation}
\eta(a \o h) = a\0 \o a\1 S(h).
\end{equation} 
Then $\beta' = \eta \circ \beta$ and $\beta: A \o_B A \to \overline{A \o H}$ is respectively injective, surjective or bijective iff
$\beta'$ is injective, surjective or bijective onto
$$ \overline{\overline{A \o H}} := \rho(1)(A \o H). $$ In particular, 
$A \| B$ is a weak
$H$-Galois extension iff $\beta':A \o_B A \to \overline{\overline{A \o H}}$ is bijective.
\end{prop}
\begin{proof}
Notice that $\overline{\overline{A \o H}}$
is a $B$-$A$-sub-bimodule of $A \o H$,
 and that $\Im \eta$ and $\Im \beta' \subseteq \overline{\overline{A \o H}}$.
Next note that an application of eq.~(\ref{eq: cdg1}) gives
\begin{eqnarray*}
\eta \beta(a \o a') & = & \eta(a{a'}\0 \o {a'}\1) \\
                    & = & a\0 {a'}\0 \o a\1 {a'}\1 S({a'}\2) \\
                    & = & a\0 {a'}\0 \o a\1 \Pi^L({a'}\1) \\
                    & = & a\0 1\0 a' \o a\1 1\1 \\
                    & = & a\0 a' \o a\1 = \beta'(a \o a').
\end{eqnarray*}
We define another linear self-mapping of $A \o H$ given
by $\overline{\eta}(a\o h) = a\0 \o \overline{S}(h)a\1$. 
Note that $\Im \overline{\eta}$ and $\Im \beta \subseteq \overline{A \o H}$.

Let $p: A \o H \to \overline{A \o H}$, $\overline{p}: A \o H \to 
\overline{\overline{A \o H}}$ be the straightforward projections given
by $p(a \o h) = a1\0 \o h1\1$, and $\overline{p}(a \o h) = 1\0 a \o 1\1 h$.
We show below that $\eta \circ p = \eta$, $\overline{\eta} \circ \overline{p} = \overline{\eta}$,
$\eta \circ \overline{\eta} = \overline{p}$ and $\overline{\eta} \circ \eta = p$, from
which it follows that the restrictions of $\eta$, $\overline{\eta}$ to $\overline{A \o H}$, $\overline{\overline{A \o H}}$
are inverses to one another, so that there is a commutative triangle connecting $\beta$, $\beta'$ via
$\eta$.  

$$\begin{diagram}
& & A \o_B A & & \\
&\SW_{\beta'} & & \SE_{\beta} & \\
\overline{\overline{A \o H}} && \lTo^{\cong}_{\eta} && \overline{A \o H} 
\end{diagram}$$

First, we note that $\eta \circ p = \eta$ since
\begin{eqnarray*}
\eta(a 1\0 \o h 1\1) & = & a\0 1\0 \o a\1 1\1 S(h 1\2) \\
                       & = & a\0 1\0 \o a\1 \Pi^L(1\1) S(h) \\
                       & = & a\0 \o a\1 S(h) = \eta(a \o h)
\end{eqnarray*}
by eqs.~(\ref{eq: bns1}) and~(\ref{eq: cdg2}).

Secondly, we note that $\overline{\eta} \circ \overline{p} = \overline{\eta}$ since
\begin{eqnarray*}
\overline{\eta}(1\0 a \o 1\1 h) & = & 1\0 a\0 \o \overline{S}(h) \overline{S}(1\2)1\1 a\1 \\
                                & = & 1\0 a\0 \o \overline{S}(h) \overline{\Pi}^R(1\1) a\1 \\
                                & = & 1\0 a\0 \o \overline{S}(h) 1\1 a\1 = \overline{\eta}(a \o h)
\end{eqnarray*}
by eqs.~(\ref{eq: over1}) and~(\ref{eq: cdg3}).  

Next we note that $\overline{\eta} \circ \eta = p$ since
\begin{eqnarray*}
\overline{\eta}(a\0 \o a\1 S(h)) & = & a\0 \o \overline{S}(a\2 S(h)) a\1 \\
                                 & = & a\0 \o h \overline{\Pi}^R(a\1) \\
                                 & = & a 1\0 \o h 1\1 = p(a \o h)
\end{eqnarray*}
by eqs.~(\ref{eq: over1}) and~(\ref{eq: cdg3}).  

Finally we note $\eta \circ \overline{\eta} = \overline{p}$ since
\begin{eqnarray*}
\eta(a\0 \o \overline{S}(h) a\1) & = &  a\0 \o a\1 S(a\2) h \\
                                & = & a\0 \o \Pi^L(a\1) h = \overline{p}(a\o h)
\end{eqnarray*}
by eq.~(\ref{eq: cdg1}). 
\end{proof}

Again let $H$ be a finite dimensional weak Hopf algebra. Recall that the $k$-dual $H^*$ is also a weak Hopf algebra,
by the self-duality of the axioms, and acts on $H$ by the usual right action $x \leftharpoonup \psi =$ $ \psi(x\1) x\2$
and a similarly defined left action. In addition, a right $H$-comodule algebra $A$ corresponds
to a left $H^*$-module algebra $A$ via $\psi \cdot a := a\0 \psi(a\1)$ \cite{NV}.   
Following Kreimer-Takeuchi and Schneider, there are two proofs that surjectivity of $\beta$ is all
that is needed in the definition of a weak Hopf $H$-Galois extension \cite{Bo, BTW}.  As a corollary of the proposition,
we offer a third and direct proof. 

\begin{cor} \cite{Bo, BTW}
\label{cor-surj}
Let $A$ be a right $H^*$-comodule algebra and $B$ its subalgebra of coinvariants $A^{\rm co \, H^*}$.
If $\beta: A \o_B A \to \overline{A \o H^*}$ is surjective, then the natural module $A_B$ is f.g.\
projective and $A | B$ is a weak $H^*$-Galois extension.
\end{cor}
\begin{proof}
We know from \cite{V} that $H$ and $H^*$ are both Frobenius algebras with nondegenerate left integral 
$t \in H$ satisfying $ht = \Pi^L(h)t$ for all $h \in H$
as well as $t \leftharpoonup T = 1_H$ for some $T \in H^*$.  Since 
$\beta$ is surjective, there are finitely many paired elements $a_i, b_i \in A$ such that
$$ 1\0 \o T 1\1 = \sum_i a_i {b_i}\0 \o {b_i}\1. $$
Let $\phi_i(a) := t \cdot (b_i a)$ for every $a \in A$.  
Then $\{ a_i \}$, $\phi_i$ are dual bases for the module $A_B$ by a computation that $\sum_i a_i \phi_i(a) = a$
for all $a \in A$, 
almost identical with \cite[p.\ 132]{Mo} for Hopf algebra actions (using the identity
$1\0 a\0 \o 1\1 a\1 = a\0 \o a\1$ at one point).

Finally, one shows that $\beta'$ is injective, for if $\sum_j u_j \o v_j \in \ker \beta'$, we
compute $$\sum_j u_j \o v_j = \sum_{i,j} a_i \o \phi_i(u_j)v_j = \sum_{i,j} a_i \o (t\1 \cdot b_i) {u_j}\0 v_j \bra {u_j}\1 , t\2 \ket = 0 $$
as in \cite[p.\ 132]{Mo}. By the proposition, $\beta$ is then
injective, whence a bijection of $A \o_B A$ onto $\overline{A \o H}$. 
\end{proof}

We next offer an example of weak Hopf-Galois extension with an
alternative proof.  For example, if $H = M_n(k)$ considered above,
the Galois map $\beta = $ \newline
$(\mu \o \id)\circ (\id \o \cop)$ given by  $\beta(e_{ij} \o e_{jk}) =$ $e_{ik} \o e_{jk}$
with coinvariants $H^L$ the diagonal matrices and $1\1 \o 1\2 =$
$\sum_i e_{ii} \o e_{ii}$, is an isomorphism by a dimension count. The
general picture is the following: 

\begin{cor} \cite[2.7]{CDG}
\label{cor - olary}
Define a coaction on $H$ by $a\0 \o a\1 = a\1 \o a\2$
for $a \in H$.
Then $H$ is a weak Hopf $H$-Galois extension of its coinvariants
$H^L$.
\end{cor}
\begin{proof}
We note that $H^L \subseteq H^{\rm co \, H}$ since
$\cop(x^L) = 1\1 x^L \o 1\2$ for $x^L \in H^L$ \cite[2.7a]{BNS}.  The converse follows
from $x \in H^{\rm co \, H}$ implies $$ x = \eps(x\1) x\2
= \eps(1\0 x) 1\1 \in H^L.$$

Next we note that $\beta'$ factors into isomorphisms in the following
commutative diagram:

$$\begin{diagram}
H \o_{H^L} H &&  \rTo^{\beta'}&& \overline{\overline{H \o H}} \\
\dTo_{\cong}^{q} && && \uTo^{\cong}_{\tau \circ (S \o S)}    \\
\overline{\overline{H \o H}} &&  \rTo_{\overline{\eta}}^{\cong}   && \overline{H \o H}
\end{diagram}$$

where $q(x \o y) :=$ $ \overline{p}(\overline{S}(x) \o_k y)$
is well-defined since $S(1\1) \o 1\2$ is a separability element
for the separable $k$-algebra $H^L$ \cite[prop.\ 2.11]{BNS}.  Its inverse
is given by $q^{-1}(\overline{p}(x \o y))$ $= S(x) \o y$. 
The mapping $\tau \circ (S \o S)$ has an obvious inverse
and is well-defined since $S(1) = 1$ and $S$ is an anti-coalgebra
homomorphism.  
\end{proof}

We provide the complete proof that a weak Hopf $H$-Galois extension
is depth two \cite[3.2]{LK}:
\begin{cor}
A weak $H$-Galois extension $A \| B$ is right and left depth two.
\end{cor}
\begin{proof}
The algebra extension $A \| B$ is right D2 since the Galois mapping $\beta: A \o_B A \stackrel{\cong}{\longrightarrow}$
$\overline{A \o H}$ and the projection $p: A \o H \to \overline{A \o H}$ are $A$-$B$-bimodule morphisms \cite[3.1]{Karl}.
Whence $A \o_B A$ is $A$-$B$-isomorphic to a direct summand $\Im p$ within $\oplus^n A$ where $n = \dim H$.

Similarly $A \| B$ is left D2 since the alternate Galois isomorphism $\beta'$  and projection $\overline{p}$
onto $\overline{\overline{A \o H}}$ 
are both $B$-$A$-bimodule morphisms.
\end{proof}

The proof of the corollary sidesteps the problem of showing $A \| B$ is a Frobenius extension, another interesting problem
for someone else or another occasion,
and implying left D2 $\Leftrightarrow$ right D2.  It is likely that a weak $H$-Galois extension is Frobenius
since there are nondegenerate integrals in $H^*$ which would define a Frobenius homomorphism via the dual
action of $H^*$ on $A$ (with invariants $B$). 
In addition we have avoided starting only with a weak bialgebra
having Galois action on $A$ and showing  the existence of an antipode on $H$ 
 in extension of \cite{Sch} for Galois actions of bialgebras. 
If we denote 
\begin{equation}
\label{eq: notation}
\sum_i \ell_i(h) \o_B r_i(h) := \beta^{-1}(1\0 \o h 1\1),
\end{equation}
we note that
\begin{eqnarray}
1\0 \o 1\1 S(h) & = & \eta(1 \o h) \\
                & = & \beta' (\beta^{-1}(1\0 \o h 1\1)) \\
                & = & \sum_i {\ell_i(h)}\0 r_i(h) \o {\ell_i(h)}\1,
\end{eqnarray}
which can conceivably be made to descend to a formula for the antipode of $H$ in terms
of just the isomorphism $\beta$.  

We  then propose  two problems and provide some evidence for each.

\begin{prob}
If $H$ is a weak Hopf algebra and $A \| B$ is $H$-Galois, is
$A \| B$  a Frobenius extension (of the second kind \cite{K})?
\end{prob}

For example, if $H$ is a Galois extension of $H^L$ as in corollary~\ref{cor - olary}
we expect such a Frobenius extension based on Pareigis's theorem that a Frobenius subalgebra $B$ of a Frobenius algebra $A$,
where the natural module $A_B$ is f.g.\ projective and the Nakayama automorphism of $A$ stabilizes $B$,
yields a $\beta$-Frobenius extension $A \| B$ where $\beta$ is the relative Nakayama automorphism of $A$ and $B$ (restrict
one and compose with the inverse of the other).  For example, if $H$ has an $S$-invariant nondegenerate integral,
the Nakayama automorphism is $S^2$ \cite[3.20]{BNS}, also the Nakayama automorphism of $H^L$, so $\beta = \id$
and $H \| H^L$ is an ordinary Frobenius extension.  

\begin{prob}
\label{conj: b} 
If $H$ is a weak bialgebra and $A \| B$ is $H$-Galois, is 
$H$ necessarily a weak Hopf algebra? 
\end{prob}

Again this is true in the special case of the weak Hopf-Galois extension in corollary~\ref{cor - olary},
a result in \cite[Brzezinski-Wisbauer]{BW}; we give another proof which may extend to  the general problem.
Note that the definition of weak Hopf-Galois extension does not make use of an antipode
nor does $H^{\rm co \, H} = H^L$ in corollary~\ref{cor - olary}.

\begin{theorem} \cite[36.14]{BW}
Let $H$ be a weak bialgebra. If the right $H$-coalgebra
$H$ with coaction $\varrho = \cop_H$ is Galois over $H^L$, then
$H$ is a weak Hopf algebra.
\end{theorem}

\begin{proof}
In terms of the notation in eq.~(\ref{eq: notation}) we define an antipode $S: H \to H$ by
\begin{equation}
S(h) = \sum_i \eps({\ell_i(h)}\1 r_i(h)) {\ell_i(h)}\2
\end{equation}

Note that by eq.~(\ref{eq: pi-ell}), $S(h) =$ $ \sum_i \ell_i(h) \Pi^L(r_i(h))$ for $h\in H$.
In order to prove that $S$ satisfies the three eqs.~(\ref{eq: bns2}),~(\ref{eq: bns1}) and~(\ref{eq: bnsa}),
we  note the three equations below for a general right $H$-comodule algebra $A$ over a weak bialgebra $H$
where $A$ is $H$-Galois over its coinvariants $B$;
the proofs are quite similar to those in \cite{Sch}.
\begin{eqnarray}
\sum_i \ell_i(h) \o {r_i(h)}\0 \o {r_i(h)}\1 & = & \sum_i \ell_i(h\1) \o r_i(h\1) \o h\2 \label{eq: s1}  \\
\sum_i a\0 \ell_i(a\1) \o_B r_i(a\1) & = & 1 \o_B a  \label{eq: s2} \\
\sum_i \ell_i(h) r_i(h) & =  & 1\0 \eps(h 1\1 ) \ \ \ \ (\forall \, a \in A, h \in H) \label{eq: s3}
\end{eqnarray}

Next we note three equations in $A \o H$, two of which we need here (and  all three might play a role in an answer to problem~\ref{conj: b}). 
\begin{eqnarray}
\sum_i {\ell_i(h\1)}\0 r_i(h\1 ) \o {\ell_i (h\1 )}\1 h\2 & = & 1\0 \o 1\1 \Pi^R(h) \label{eq: ex1} \\
\sum_i {\ell_i(h\2)}\0 r_i(h\2) \o h\1 {\ell_i(h\2 )}\1 & = & 1\0 \o \Pi^L(h 1\1 ) \label{eq: ex2}
\end{eqnarray}

\begin{equation}
\label{eq: ex3}
\sum_i {\ell_i(h\1 )}\0 r_i(h\1 ) {\ell_i(h\3)}\0 r_i (h\3 ) \o {\ell_i(h\1 )}\1 h\2 {\ell_i(h\3 )}\1 = 
\end{equation} 
$$ \sum_i {\ell_i(h)}\0 r_i(h) \o {\ell_i(h)}\1. $$ 
They are established somewhat similarly to \cite{Sch} and left as exercises.

Applying eq.~(\ref{eq: ex1}) with $A = H$ and $a\0 \o a\1 = a\1 \o a\2$, we obtain one of the antipode axioms:
\begin{eqnarray*}
S(h\1) h\2 & = & \sum_i \eps({\ell_i(h\1)}\1 r_i(h\1 )){\ell_i (h\1 )}\2 h\2  \\
           & = & \eps(1\1 ) 1\2 \Pi^R(h) = \Pi^R(h). \ \ \ \ (\forall \, h \in H) 
\end{eqnarray*}
 
Applying eq.~(\ref{eq: ex2}), we obtain
\begin{eqnarray*}
h\1 S(h\2) & = & \sum_i \eps({\ell_i(h\2)}\1 r_i(h\2) )h\1 {\ell_i(h\2 )}\2 \\
           & = & \eps(1\1 ) \Pi^L(h 1\2 ) = \Pi^L(h) \ \ \ \ (\forall \, h \in H) 
\end{eqnarray*}

Finally we see $S$ is an antipode from the just established eq.~(\ref{eq: bns2}) and applying eq.~(\ref{eq: pi-are}): 
\begin{eqnarray*}
\Pi^R(h\1) S(h\2) & = & \sum_i \Pi^R(h\1) \ell_i(h\2) \Pi^L(r_i(h\2)) \\
                   &=&  \sum_i 1\1 \ell_i(h1\2) \Pi^L(r_i(h1\2 )) \\
                   & = & \sum_i \ell_i(h) \Pi^L(r_i(h)) = S(h)  
\end{eqnarray*}
where we use the general fact that $\beta$ is left $A$-linear, so $\sum_i 1\0 \ell_i(h 1\1) \o r_i(h 1\1) = $
$\beta^{-1}({1'}\0 1\0 \o h {1'}\1 1\1) = $ $ \sum_i \ell_i(h) \o r_i(h)$.  
\end{proof}

 


\begin{thebibliography}{XXXXXX}
\begin{small}
\bibitem{BaS}{I. B\'{a}lint and K.~Szlach\'anyi,
Finitary Galois extensions over non-commutative bases,
in preparation.}
\bibitem{BTW}{T.~Brzezi\'nski, R.B.~Turner and A.P.~Wrightson,
The structure of weak coalgebra-Galois extensions,
preprint (2004), \texttt{QA/0411230}.}
\bibitem{BW}{T.~Brzezi\'nski and R.~Wisbauer,
\textit{Corings and Comodules}, LMS \textbf{309}, Cambridge University Press, 2003.}
\bibitem{Bo}{G.~B\"{o}hm,
Galois theory for Hopf algebroids, KFKI preprint (2004), \texttt{RA/0409513}.}
\bibitem{BS}{G.~B\"ohm and K.~Szlach\'anyi,
Hopf algebroids with bijective antipodes: axioms, integrals and duals,
\textit{J. Algebra} \textbf{274} (2004), 708--750.}
\bibitem{BNS}{G.~B\"ohm, F.~Nill and K.~Szlach\'anyi,
Weak Hopf algebras, I. Integral theory and $C^*$-structure,
\textit{J. Algebra} \textbf{221} (1999), 385-438.}
\bibitem{CDG}{S.~Caenepeel and E. De Groot,
Galois theory for weak Hopf algebras, preprint (2004),
RA/0406186.}
\bibitem{EN}{P.~Etingof and D.~Nikshych,
Dynamical quantum groups at roots of $1$, 
\textit{Duke Math. J.}, \textbf{108} (2001), 135--168.}
\bibitem{FMS}{D.~Fischman, S.~Montgomery, and H.-J.~Schneider,
Frobenius extensions of subalgebras of Hopf algebras,
{\it Trans.\ Amer.\ Math.\ Soc.\ } {\bf 349} (1997), 4857--4895.}
\bibitem{K}{L.~Kadison,
\textit{New Examples of Frobenius Extensions}, University Lecture Series \textbf{14}, AMS, Providence, 1999.
Update, 6 pp: www.ams.org/bookpages.}
\bibitem{Karl}{L.~Kadison,
Depth two and the Galois coring, preprint (2004), \texttt{RA/0408155}.}
\bibitem{LK}{L.~Kadison,
An action-free characterization of weak Hopf-Galois extensions,
preprint (2004), \texttt{QA/0409589}.}
\bibitem{KK}{L. Kadison and B.~K\"ulshammer,
Depth two, normality and a trace ideal condition
for Frobenius extensions,  \texttt{GR/0409346}.}
\bibitem{KN}{L.~Kadison and D.~Nikshych,
Hopf algebra actions on strongly separable extensions of depth two,
\textit{Adv.\ in Math.} \textbf{163} (2001), 258--286.}
\bibitem{KSt} L.~Kadison and A.A.~Stolin, 
             Separability and Hopf algebras,
             in: \textit{Algebra and its Applications},  
eds.\ D.V.~Huynh, S.K.~Jain, and S.~Lopez-Permouth, Contemporary  Math.\ \textbf{259}     
    A.M.S., Providence, 2000, 279-298. 
\bibitem{KS}{L.~Kadison and K.~Szlach\'anyi,
Bialgebroid actions on depth two extensions and duality, \textit{Adv.\
in Math.} \textbf{179} (2003), 75--121. 
\texttt{RA/0108067}}
\bibitem{Lu}{J.-H. Lu,
Hopf algebroids and quantum groupoids,
\textit{Int.\ J. Math.} \textbf{7} (1996), 47--70.}
\bibitem{Maj}{S.~Majid,
\textit{Foundations of Quantum Group Theory}, Cambridge Univ.\ Press, 1995.}
\bibitem{Mo}{S.~Montgomery,
\textit{Hopf Algebras and Their Actions on Rings},  CBMS Regional Conf.\ 
Series in Math.\ 
Vol.\ 82, AMS, Providence, 1993.}
\bibitem{NV} D.~Nikshych and L.~Vainerman,
Finite dimensional quantum groupoids
and their applications, in:  \textit{New Directions in Hopf Algebras},
S.~Montgomery and H.-J.~Schneider, eds., MSRI Publications, vol.\ \textbf{43}, Cambridge, 2002, pp.\ 211--262. 
\bibitem{Sch}P.~Schauenburg,
A bialgebra that admits a Hopf-Galois extension is a Hopf algebra,
Proc.\ A.M.S.\ 125 (1997), 83--85. 
\bibitem{HJS}{H.-J. Schneider,
Normal basis and transitivity of crossed products for Hopf algebras,
{\it J. Algebra} {\bf 151} (1992), 289--312.}
\bibitem{Sz}{K.~Szlachanyi, Galois actions by finite quantum groupoids,
\textit{Locally compact quantum groups and groupoids (Strasbourg, 2002)},IRMA Lect.\ Math.\ Phys.\ 2, de Gruyter, Berlin, 2003, 105--125. \texttt{QA/0205229}. }
\bibitem{V}{P.~Vecserny\'es, Larson-Sweedler theorem and the role of grouplike elements in weak Hopf algebras,
\texttt{QA/0111045}.} 

\end{small}
\end{thebibliography}
\end{document}